# THE SPACE OF PLANAR SOAP BUBBLE CLUSTERS

Frank Morgan
Williams College
Frank.Morgan@williams.edu

Soap bubbles and foams have been extensively studied by scientists, engineers, and mathematicians as models for organisms and materials, with applications ranging from extinguishing fires to mining to baking bread; see for example Cantat *et al.* [CC], Weaire-Hutzler [WH], Morgan [M, Chap. 13], and references therein. Here we provide some basic results on the space of planar clusters of $n$ bubbles of fixed topology. We show for example that such a space of clusters with positive second variation is an $n$-dimensional manifold, although the larger space without the positive second variation assumption can have singularities. Earlier work of Moukarzel [Mou] showed how to realize a cluster as a generalized Voronoi partition, though not canonically.

## 1. Soap Bubble Clusters

**1.1. Definitions.** A planar *cluster* consists of disjoint circular arcs/line segments meeting in threes at positive angles, enclosing $n$ connected regions, with areas denoted $A_i$. We assume that the cluster is connected and has at least two regions ($n \geq 2$). A cluster is in *quasi-equilibrium* if the arcs meet at 120 degrees. For *equilibrium* we further assume that the sum of the curvatures around a path from a region to itself is 0, or equivalently around a vertex. This condition makes *pressure* well defined as the sum of the curvatures along a path from the exterior to the point. An equivalent form of the condition is that circles/lines through a vertex meet again (each forward or backwards, at infinity in the case of three straight lines), as follows by the law of sines [M, Fig. 14.1.2]. All of these definitions are preserved by Möbius (linear fractional) transformations. A cluster is in equilibrium if and only if the cluster has vanishing first variation of length under smooth deformations of the plane which preserve the areas (see [Q, Prop. 2.6 and Appendix A]). A cluster is *minimizing* if it minimizes length for given areas. A minimizing cluster has vanishing first variation and nonnegative second variation.

More generally and technically one might define a soap bubble cluster in $\mathbf{R}^N$ as $n$ disjoint regions of finite volume and perimeter such that Lipschitz deformations inside small balls preserving the volumes cannot reduce the area of the union of the boundaries (see [M, §11.3]). In $\mathbf{R}^2$ such are equilibrium clusters as defined above (see [M, §13.10]).

Parameterize clusters of given combinatorial type by $v$ vertices $V_i$ and $e$ areas $a_{ij} = -a_{ji}$ between edges and straight lines between vertices. Since $e = 1.5v$, by Euler there are $2(n-1)$ vertices and $3(n-1)$ edges. On this $(7n-7)$-dimensional manifold $C$ the areas $A_i$ and perimeter $P$ are smooth functions. (Using areas $a_{ij}$ instead of curvatures avoids ambiguity of large and small circular arcs of the same curvature.) Since rigid motions infinitesimally have no fixed points (because by the assumption $n \geq 2$ we exclude a circle ($n=1$)), there is a smooth quotient manifold $Q$ of dimension $7n-10$.

**1.2. Lemma.** *The Jacobian of the area vector $(A_1, A_2, \ldots, A_n)$ has full rank n (on C and on Q).*

*Proof.* The area of a region sharing an edge $E_{ij}$ with the exterior can be varied by varying $a_{ij}$. Adjacent regions can be similarly adjusted by varying the $a_{ij}$ of a shared edge and restoring the area of the outer region. Work your way inward to obtain arbitrary variations.

**1.3. Corollary**. *Level sets of fixed area vector in the spaces of clusters (C and Q) provide a smooth foliation (by the Inverse Function Theorem).*

**1.4. Theorem.** *The space in Q of equilibrium clusters modulo rigid motions with n regions with positive second variation for fixed areas is a smooth n-dimensional manifold, locally parametrized by the areas.*

*Proof.* Varying areas smoothly preserves positive second variation equilibria.

**1.5. Conjecture.** *An equilibrium of positive second variation is unique for given area vector and given combinatorial type.*

**1.6. Remarks.** The theorem and conjecture do not hold with positive second variation replaced by nonnegative second variation or stability, as observed by Weaire *et al.* [WCG]. As a trivial example, consider a bubble with two small lenses in its boundary, with the extra degree of freedom of the distance between the lenses, yielding a (3+1)-dimensional manifold of stable equilibria. Lenses tend to add an extra dimension. In the symmetric case, the two middle arcs can be replaced by arcs of a different curvature, yielding quasi-equilibria, as pointed out to me by John M. Sullivan, another 4-dimensional manifold containing the symmetric case, which is thus a non-manifold point in the space of quasi-equilibria.

For a more interesting example, a long chain or "necklace" of $n-1$ say unit curvature bubbles surrounding a chamber with 0 pressure as in Figure 1 is floppy, with a multi-parameter family of configurations with the same areas and pressures. To compute the dimension, note that you can slide each bubble around the adjacent one, for $n-2$ free parameters, minus 2 for them to match up and 1 for rigid rotation and 1 more to preserve the area of the central chamber, for a total of $n-5$. These necklaces can probably be shown minimizing by the methods of Cox *et al*. [CH]. The symmetric one also sits in a one-parameter family of varying area and pressure of the central chamber, apparently not available in the non-symmetric case, so that the space is not a manifold at that point. For smaller area the pressure goes negative and the cluster becomes unstable. For larger area the pressure goes positive, the cluster becomes strictly stable, and all of the areas can be varied, but the cluster is no longer floppy.

Is the set of equilibria of fixed combinatorial type connected? Equilibria with nonnegative second variation? For fixed area vector?

The theorem and proof hold if length is replaced by any smooth, uniformly convex norm. For equilibrium, curves meet in threes such that the unit duals of the tangent vectors and the (constant) generalized curvatures both sum to zero. See [Le], [MFG].

It is unknown whether the manifold of Theorem 1.4 is locally parametrized by pressures, except for double bubbles in $\mathbf{R}^2$, where the only stable double bubble is the standard double bubble [MW].

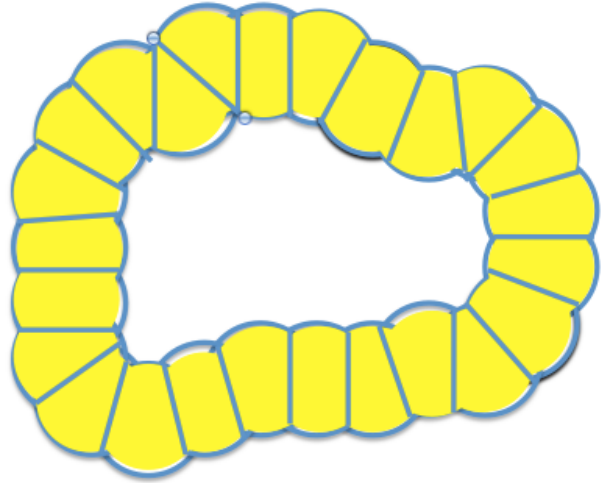

Figure 1. A long chain or "necklace" of $n-1$ say unit curvature bubbles surrounding a chamber with 0 pressure is floppy.

## 2. Dimension of Equilibria

An upper bound on the dimension on the space of equilibria is given by placing the vertices ($2v$) and the orientation of the curves leaving one vertex (1), which determines the orientations at adjacent vertices, etc., minus rigid motions (3), for a total of $2v-2 = 4n-6$. Fenyes [F] proved that if every bubble has at least three sides the "$s_{ij}$" conditions on each edge that the angle of the curve with the chord is the same at both ends are independent on vertices plus orientations, which yields a bound of $3v-3 -1.5v = 1.5v-3 = 3n-6$. Adding $v-1$ curvature constraints would bring the estimate down to $.5v-2 = n-3$, just 3 less than expected, probably due to some minor dependence among the curvature constraints. Indeed, none are needed for the triple bubble ($v=4$), and using $v-4$ instead of $v-1$ would yield exactly the expected $n$. So probably generically the dimension is due to placing $v$ vertices, choosing directions at each vertex, the $s_{ij}$ edge constraints, the curvature constraints, and modding out by rigid motions:

$$2v + v - 1.5v - (v-4) - 3 = .5v + 1 = n.$$

To prove Fenyes's lemma, sending the conjugate vertex to infinity, one has straight lines from points $x_i$, $x_j$, and $x_k$ meeting $x_h$ at 120 degrees. The result holds because for fixed $x_i$, $x_j$, $x_k$ and directions from each of them, the point $x_h$ is the only point where circular arcs from $x_i$, $x_j$, and $x_k$ in those directions can meet at 120 degrees.

*Question.* Is the space of equilibria of given combinatorial type connected? for given area vector?

# 3. Small Clusters

Any two vertices support a 1-parameter family of quasi-equilibrium double bubbles, all equilibria by trigonometry, all minimizing by the planar double bubble theorem (see [M]). The number of parameters equals 4 – 3 [rigid motions] = 1, and 1+1 = 2.

Any three vertices uniquely determine the fourth vertex and a triple quasi-equilibrium bubble (easy in the most symmetric case, in general by Möbius transformations), which in fact are equilibria. By Wichiramala [Wi], they are minimizing. The number of parameters is 6 – 3 = 3 (or the previous 2 plus the 1-parameter family of decorations).

Similarly any such three vertices and an appropriate fourth on the triple bubble determine a standard 4-bubble. The number of parameters equals 6 + 1 – 3 = 4 (or the previous 3 plus 1).

A 4D family of equilibrium 3-clusters is given by a circle decorated by two lenses. In the symmetric case, the two middle arcs can be replaced by arcs of a different curvature, as pointed out to me by John M. Sullivan, yielding another 4D family of quasi-equilibria.

The standard type 4-cluster similarly admits quasi-equilibria. Start with an equal-area double bubble, with one bubble above the other. Decorate the two vertices symmetrically to produce a 4-bubble symmetric under horizontal as well as vertical reflection. Now lengthen the horizontal line in the middle and replace the top and bottom arcs with others of smaller curvature to maintain the 120-degree angles.

The first equilibrium cluster not to arise from repeated decorations is the flower with four petals around a 4-sided central bubble as in Figure 2 below. Is it part of a 5-dimensional family of equilibria? (Yes, next paragraphs.) Are they stable? (I think so.) Are there more quasi-equilibria? (Yes, similar to previous constructions.)

As Sullivan (email 2014) pointed out, "it's easy to find a 4-parameter family: note that in the symmetric cluster, the size of the inner bubble can be varied; then apply Möbius transformations including scaling. What this doesn't give is clusters of the form where the four bubbles in the ring around the center alternate large, small, large, small [because any inversion keeping one pair of opposites equal makes the other pair unequal]. But I think those could be constructed easily explicitly (with $D_2$ symmetry)."

An easy way to construct these 5-clusters with $D_2$ symmetry is to note that a quarter is part of a double bubble with a lens (which may extend outside the double bubble). So start with a double bubble, with one bubble above the other, with a centered lens (three parameters) and slide the lens until the line perpendicular to two upper circular arcs is perpendicular to the line perpendicular to the two lower circular arcs. Now reflection across the two lines yields the desired 5-cluster. At least for nearly equal petal areas, these clusters are all stable [B]. Now Möbius inversion yields a 5-parameter family of equilibrium 5-clusters (without imposed symmetry).

The space of circular arc triangles with three 120-degree angles has dimension 3 (vertices arbitrary by linear fractional transformations and uniqueness obvious from equilateral case).

# 4. 3D Clusters

The space of clusters in space is harder to parameterize; we follow the methods of White [Wh] (see his Intro. and Thm. 3.1). We consider clusters defined as smooth surfaces meeting in threes at positive angles along smooth curves, which in turn meet in fours at positive angles, enclosing $n$ connected regions, with volumes denoted $V_i$. An equilibrium cluster consists of constant-mean curvature surfaces (not necessarily spherical) meeting in threes at 120 degrees along smooth curves, which in turn meet in fours at equal angles of $\cos^{-1}(-1/3) \approx 109$ degrees (see [M, 11.3]). We assume that a cluster is connected and has at least two regions ($n \geq 2$). At a given equilibrium cluster $C$, we consider perturbations as follows. First we perturb the singular points, replacing the singular curves by homothetic copies. Next via a real-analytic transverse vectorfield we identify a tubular neighborhood of each singular curve of length $L$ with $[0, L] \times \mathbf{B}(0, \varepsilon)$, and perturb it by a smooth function from $L$ to $\mathbf{B}(0, \varepsilon)$ vanishing at the endpoints. Having thus perturbed the singular curves, we adjust each component of the surface by the harmonic function with the given change in boundary values. Finally, we perturb the surface by a smooth normal variation vanishing at the boundary. This process parameterizes a neighborhood in the space of clusters of a given cluster with a neighborhood of the origin in a space of smooth functions.

We take as our space of smooth functions the Banach space $C^{2,\alpha}$ (fixed $0 < \alpha < 1$) of twice Hölder differentiable functions under the standard $C^{2,\alpha}$ norm, making the space of clusters a smooth Banach manifold. The enclosed volume vector is a real-analytic function on this linear space. The advantage of the Hölder spaces is that harmonic functions inherit the smoothness of the boundary values up to the boundary. White [Wh, §1.5] explains how to modify the space to make it separable. A good reference on Banach manifolds is provided by Lang [L].

**4.1. Theorem.** *The space of clusters of n bubbles modulo rigid motions with positive second variation in $\mathbf{R}^N$ is a smooth n-dimensional manifold, locally parametrized by the volumes.*

*Proof.* We consider perturbations as above of a fixed cluster with positive second variation. As in the 2D case, now by the inverse function theorem on Banach manifolds, the space is $C^\infty$ foliated by submanifolds of fixed volume vector. The theorem follows. (Because the range of the volume vector $V$ is finite dimensional, there is trivially a subspace complementary to ker $DV$, yielding the desired local isomorphism of the tangent space with ker $DV \times \mathbf{R}^n$, which by the inverse function theorem implies that the space of clusters is locally diffeomorphic to $\{V=0\} \times \{V\}$.)

*Remarks.* It is unknown whether the manifold is locally parametrized by pressures, even for double bubbles in $\mathbf{R}^3$, where it was a major advance to prove that pressure is strictly decreasing in volume for *area-minimizing* double bubbles ([H, Thm. 3.2], see [M, 14.5]).

The above approach and Theorem 4.1 hold for piecewise smooth clusters (stratified manifolds) in $\mathbf{R}^N$. More generally and technically one could define an (equilibrium) soap bubble cluster in $\mathbf{R}^N$ as $n$ disjoint regions of finite volume and perimeter such that Lipschitz deformations inside small balls preserving the volumes cannot reduce the area of the union of the boundaries (see [M, 11.3]). In $\mathbf{R}^2$ and $\mathbf{R}^3$ such are equilibrium clusters as previously defined, although in $\mathbf{R}^3$ it is not known whether every four-curve singularity satisfies the local area minimization condition (see [LM]). In higher dimensions, it is not known whether such general equilibrium clusters are piecewise smooth (stratified manifolds).

## 5. de Sitter Spacetime

The following theorem was communicated to me by Aaron Fenyes, who also provided Figures 2 and 3.

**5.1. Theorem** (Fenyes, 2014). *There is a natural correspondence between connected n-junction immersed equilibrium clusters in the plane with a point at infinity and the algebraic variety of n triples of points in de Sitter spacetime, each triple evenly spaced at distances $2\pi/3$ on an oriented line, the 3n (not necessarily distinct) points in antipodal pairs.*

For simplicity of description, we are assuming that the junctions of the cluster are distinct and that the triples on oriented lines are distinct (although for example you could have the same triple evenly spaced on the same line with *opposite* orientations, corresponding to three films meeting at two points, as in the double bubble). The triple associated to a junction consists of the oriented circles/lines leaving the junction, the orientation of their line giving their counterclockwise order around the junction.

Figures 2 and 3 provide a schematic representation in de Sitter spacetime of the 4-flower discussed above and a two-parameter family of deformations.

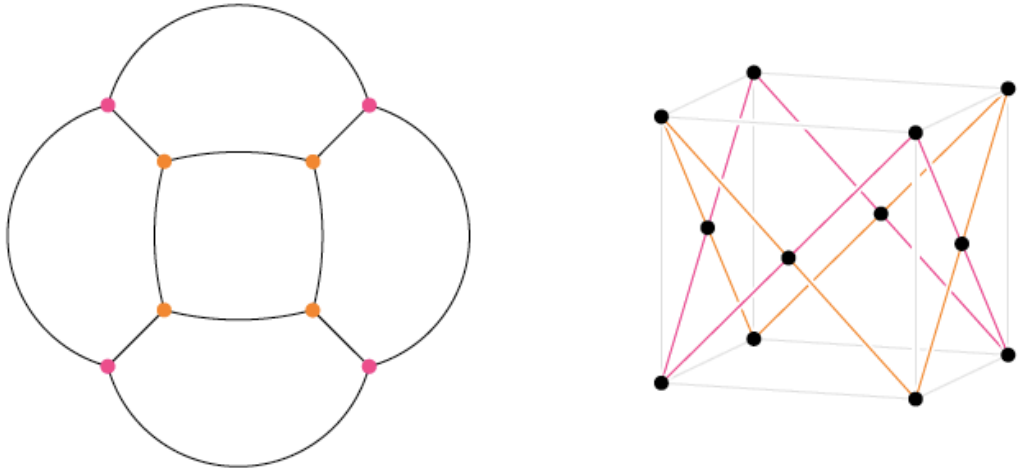

Figure 2. A 4-flower and a schematic of its representation in de Sitter spacetime, where points (representing planar circles) occur in equally spaced triples on geodesics by the equilibrium conditions at each junction of the soap film. Aaron Fenyes.

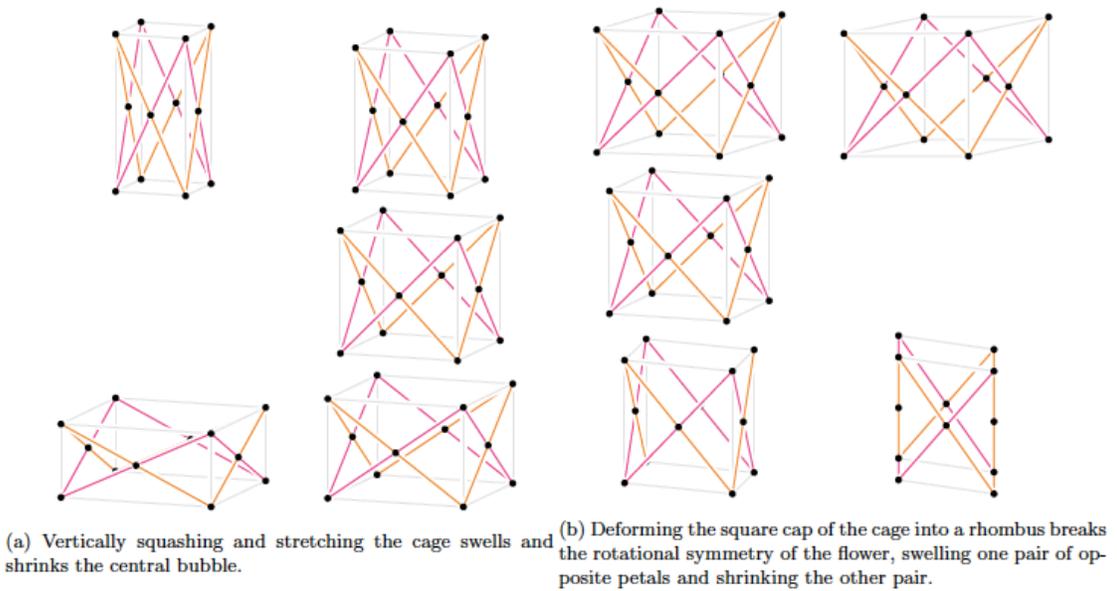

(a) Vertically squashing and stretching the cage swells and shrinks the central bubble.

(b) Deforming the square cap of the cage into a rhombus breaks the rotational symmetry of the flower, swelling one pair of opposite petals and shrinking the other pair.

Figure 3. Two families of deformations of the representation in de Sitter spacetime of the 4-flower. Aaron Fenyes.

*Proof of Theorem* 5.1. de Sitter spacetime S consists of the oriented spacelike lines in

$M$ = (3+1)-dimensional Minkowski spacetime = {2×2 Hermitian matrices}

(see [S]). Three oriented circles/lines meet at two points at $2\pi/3$ radians if and only if they are collinear and evenly spaced at distances $2\pi/3$ in de Sitter spacetime, and the orientation of the line determines their cyclic order and hence which of the two points is in the cluster. The pairing of the points in antipodal pairs corresponds to the fact that each circle/line leaves one junction and enters another; hence appears with both orientations, antipodal points in *S*.

*Acknowledgements.* These discussions began with undergraduate research by Aaron Fenyes [F] and continued when Morgan and John M. Sullivan were attending the inspiring conference "Foams and Minimal Surfaces—12 Years On" newton.ac.uk/programmes/FMS/fmsw02.shtml at the Isaac Newton Institute, 24-28 February 2014, organized by Simon Cox and Denis Weaire. I am grateful to Fenyes and Sullivan for their contributions.

# Appendix on Decoration

**A.1. Proposition.** *A 3-sided bubble in a cluster can be expanded or shrunk (even to a point) without affecting the rest of the cluster.*

*Proof sketch.* Given a 3-sided bubble, use a linear fractional transformation to map a vertex to infinity and the other point where the two incident edges meet to the origin. The third edge is mapped to one side of an equilateral 3-sided bubble centered at the origin. Now map its third vertex to infinity. Our original bubble is now mapped to an equilateral 3-sided bubble, which can be shrunk or expanded without affecting the rest of the cluster.

To show failure for *k*-sided bubbles, we need the following lemma about double bubbles. We know from double bubbles that when a circular arc splits into two others in an equilibrium cluster and they continue until they meet again, the curve that emerges is a continuation of the original circular arc.

**A.2. Lemma.** *For a double bubble consisting of bubbles of radius $r_1$ and $r_2$, increasing $r_2$ causes the distance d between their centers to decrease if $r_2 < .5r_1$ and increase if $r_2 > .5r_1$.*

*Proof.* See Morgan [M, Fig. 14.1.2]. By the law of cosines,

$$d^2 = r_1^2 + r_2^2 - r_1 r_2.$$

Therefore $2d(dd/dr_2) = 2r_2 - r_1$, which is negative if $r_2 < .5r_1$ and positive if $r_2 > .5r_1$.

**A.3. Proposition.** *In general, a bubble cannot vary pressure without affecting the rest of the cluster.*

*Proof.* Consider a bubble in contact with three other surrounding, fixed, say smaller, bubbles. As its pressure decreases and it grows, it must move farther from all three—impossible. The other three are connected if you like by chains of smaller bubbles.

**A.4. Remark**. The easier, well-known Decoration Proposition, that you can insert a three-sided bubble at a triple junction, follows from inversion to infinity of the other point where the three curves meet, which exists by double bubble trigonometry. In 3D, such an argument, based on triple bubbles, holds only if the six surfaces meeting at the tetrahedral junction are all spherical. How do you prove a general decoration theorem in 3D?